\documentclass[12 pt]{article}
\usepackage[dvips]{graphicx}
\usepackage{amssymb}

\newcommand{\be}{\begin{equation}}
\newcommand{\ee}{\end{equation}}
\newcommand{\bea}{\begin{eqnarray}}
\newcommand{\eea}{\end{eqnarray}}
\def\reff#1{(\ref{#1})}

\newtheorem{theorem}{Theorem}

\begin{document}
\bibstyle{ams}

\title{
A densest compact planar packing \\
with two sizes of discs \\
}

\author{Tom Kennedy
\\Department of Mathematics
\\University of Arizona
\\Tucson, AZ 85721
\\ email: tgk@math.arizona.edu
\bigskip
}

\maketitle

\begin{abstract} 
We consider packings of the plane using discs of 
radius $1$ and $r=0.545151\cdots$. The value of $r$ admits compact
packings in which each hole in the packing is formed by three discs 
which are tangent to each other. We prove that the largest density
possible is that of the compact packing shown in figure \ref{fig545}.

\end{abstract}

\vspace{\fill}
\hrule width2truein
\smallskip
{\baselineskip=12pt
\noindent
\copyright\ 2004 by the author. Reproduction of
this article is permitted for non-commercial purposes.
\par}

\newpage

\section{Introduction} \setcounter{equation}{0}

We consider the following packing question in two dimensions. Fix a 
number $r<1$. Using discs of radius $1$ and $r$, 
what is the densest packing of the plane? We do not impose any constraint
on the relative number of discs of the two sizes. 

It was proved long ago that the densest packing of the plane using discs
of equal radii is to put the centers of the discs on a triangular lattice.
\cite{thua,thub}. (The Voronoi cells of this packing are hexagons, and this
packing is often referred to as hexagonal or honeycomb.)
The density of this triangular packing is $\pi/\sqrt{12}$. 
L. Fejes T\'oth observed that if $r$ is slightly less than $1$, then 
one cannot do any better than this packing density \cite{ftc}. 
The interval in which it has been proved that the highest packing density is 
$\pi/\sqrt{12}$ was increased to $[0.906\cdots,1]$ 
by  Florian \cite{florb}. Using an an idea of Boroczky [Bo], 
Blind \cite{blia,blib} and G. Fejes T\'oth \cite{gtoth} independently 
extended it to $[0.742 \cdots,1]$. 

For smaller values of $r$ there is a rich variety of packings with densities
greater than $\pi/\sqrt{12}$. 
(A survey of the best known 
packings as a function of $r$ may be found at \cite{kenb}.)
However, the densest packing has been 
rigorously established only for six particular values of $r$.
All six of these values of $r$ allow compact packings. 
A packing is said to be compact if each disc is surrounded by a ring 
of discs, all of which are tangent to the disc at the center. Furthermore,
each disc in the ring is tangent to the two discs in the ring which are 
adjacent to it in the cyclic order.  
Heppes has proved that for six values of $r$ which allow compact packings, 
the largest density is attained by a particular compact packing.
It has been shown that there are only nine values of $r$ 
that admit compact packings \cite {kena}. 

\begin{figure}[tbh]
\includegraphics{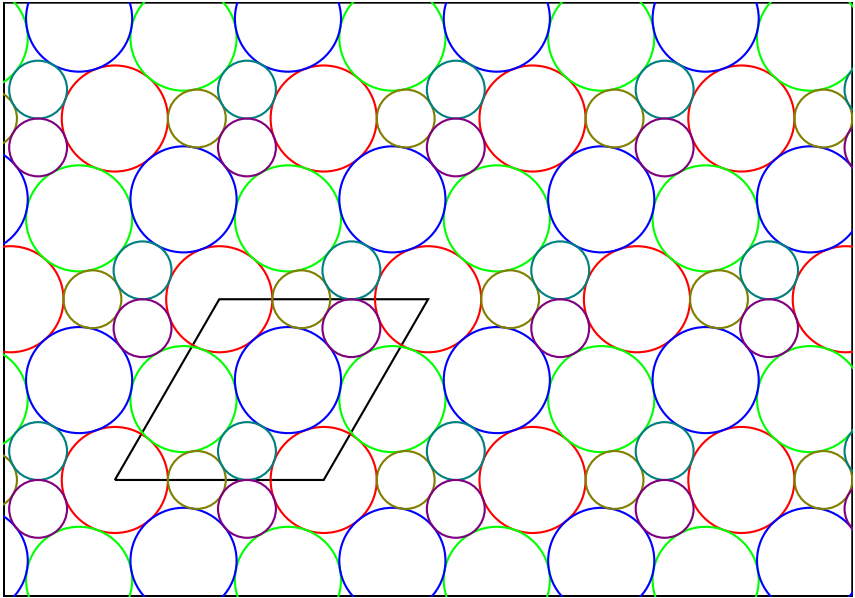}
\caption{A densest packing with $r=0.545 \cdots$. The quadrilateral shows
a unit cell for the packing.}
\label{fig545}
\end{figure}

In this paper we consider one of the values of $r$ which admits compact
packings but for which it has not been shown that 
a compact packing attains the largest density.
The value is $r=0.545151042 \cdots$. The exact
$r$ is a root of  
\be
r^8 - 8r^7 - 44 r^6 -232 r^5 -482 r^4
-24 r^3+388 r^2 -120 r + 9 = 0 
\ee
(This equation is derived in appendix \ref{appen_eqr}.)
In the remainder of this paper $r$ will denote this particular radius.
In this paper we prove that the largest packing density possible 
using discs of radius $1$ and $r$ is that attained by the compact packing 
in figure \ref{fig545}.
It has a packing density of $\delta= 0.911627478 \cdots$.   

The strategy of our proof comes from a technique in classical statistical
mechanics known as ``m-potentials'' \cite{sla}. Our strategy is similar to 
Heppes ``cell balancing'' \cite{hepb}.
In statistical mechanics m-potentials were introduced to 
deal with frustrated spin systems. In spin systems the Hamiltonian 
or energy function is typically a sum over translates of a local energy
function. A ground state of such a system is a configuration of spins 
which minimizes the total energy function. One can also ask what is 
the minimum of a single local energy function. This may be less than
the value of the local energy function in the ground state. When 
this happens the system is said to be frustrated. What one would like 
to do locally to minimize the energy cannot be done globally to 
simultaneously minimize all the local energy functions. 

In classical spin systems an m-potential is a local function 
on the spin configurations with the property that when it is summed over
all translates the result is just the zero function. Thus one can add 
this local function to the local energy function and obtain a total energy 
function with the same ground states. In some problems by carefully 
chosing the m-potential one can obtain a system which is not frustrated.
(We should note that an important class of frustrated spin systems 
comes from disordered systems such as spin glasses, but m-potentials
have not been useful in this context) 

The disc packing problem is also frustrated. Given just three discs the 
densest packing is to have them all touch one another. 
For most values of $r$ it is not possible to find a packing
in which the packing locally always consists of three discs 
touching each other. Even when it is possible (as it is for our value 
of $r$), different choices for the three discs that are tangent 
will give different densities for the triangle formed, and it is 
not possible to find a packing that only uses the densest triangle. 
By analogy with the m-potentials, we will introduce a function on 
packings that is a sum over all triangles of a function of 
three discs. We will refer to it as a ``localizing potential'' 
since its purpose is to reduce the global problem of finding the best 
packing to a local problem involving only three discs. 
Note that for the spin problem, if we are only interested in the 
ground states then it suffices that the sum of the local m-potentials
be non-negative. Likewise, in the packing problem it suffices that 
the sum of the localizing potentials be non-negative. 

In this paper we only allow discs of two sizes. An interesting and 
presumably more complicated question is what is the largest packing 
density if we allow discs with any radius in $[r,1]$.
Another interesting question is what is the densest packing if we add the 
constraint that the ratio of the number of one type of disc to the other 
must converge to a given value in the limit of packing the entire plane. 
This question was studied non-rigorously in \cite{lh}.

In the next section we explain the method of localizing potentials
in detail. Our localizing potential is the sum of two parts, 
a ``vertex localizing potential'' and an ``edge localizing potential.'' 
In section \ref{section_vertex_lp} we define the vertex localizing 
potentials used in the proof of our result and show that the sum of 
our vertex localizing potentials is non-negative. 
We define the edge localizing potentials and show their sum is 
non-negative in section \ref{section_edge_lp}.
Sections \ref{section_local} and \ref{section_global} are devoted to 
studying the local problem involving just three discs that comes 
from adding our localizing potential to the original packing problem.

\section{ Localizing potentials}
\label{section_loc_pot}

Consider the centers of the discs in a large packing.
The Delaunay decomposition gives a triangulation 
in which the vertices of the triangles are 
the centers of the discs. We denote the triangles by $T_i$.
Let $A(T)$ be the area of triangle $T$.
Given a triangle $T$, let $\phi_0,\phi_1,\phi_2$
be the angles in $T$ and $r_i$ the radii of the discs at the 
corresponding vertices. We define 
\be
D(T)= \sum_{i=0}^2 {1 \over 2} \phi_i r_i^2 
\ee
As long as the triangle is not too ``flat'', $D(T)$ is the area of 
the intersection of the triangle with the three discs.
Even when the triangulation contains triangles for which this is not 
true, the sum $\sum_i D(T_i)$ will be (up to boundary effects) 
the total area of the discs 
in the packing since the sum of the angles $\phi$ around a vertex is 
always $2 \pi$. Up to boundary effects the packing density is 
\be
{\sum_i D(T_i) \over \sum_i A(T_i)} 
\label{global_prob}
\ee
So we want to prove this ratio is no greater than $\delta$. 
This is equivalent to 
\be
\sum_i (\delta A(T_i) - D(T_i)) \ge 0
\ee
We define 
\be 
E(T) = \delta A(T) - D(T) 
\ee
Heppes defines the ``surplus area'' to be  $A(T) - D(T)/\delta$.
We refer to the quantity $E(T)$ as the ``excess'' of the triangle.
If it were nonnegative for every triangle we would be done. 
This happens to be true when one considers packings with discs of a single
radius, but it is not true in our problem. In an optimal packing,
triangles have both positive and negative excess, but the sum over all 
triangles of the excess is zero. We must prove that for any packing the 
sum over the triangles of the excess is non-negative.

We want to define a function $F(T)$ on triangles with the following 
two properties. First, for any Delaunay decomposition we require 
\be
\sum_i F(T_i) \ge 0
\label{global_lp_cond}
\ee
Second, for any triangle that can occur in a Delaunay decomposition, we require
\be
E(T) - F(T) \ge 0
\label{local_prob}
\ee
If we can do this, then we are done :
\be
\sum_i E(T_i) \ge \sum_i (E(T_i) - F(T_i)) \ge 0
\ee
We refer to $F$ as a {\it localizing potential} since it reduces 
proving the global inequality \reff{global_prob} to proving the local one 
\reff{local_prob}. We will prove the following theorem.

\begin{theorem} For the particular value of $r$ we are considering 
($r \approx 0.545151$), there is a localizing potential $F(T)$ which 
satisfies inequalities \reff{global_lp_cond} and \reff{local_prob}
with $\delta$ equal to the density of the packing shown in figure 
\ref{fig545}. Thus the density of a packing consisting of discs of 
radius $1$ and $r$ is at most $\delta$. 
\end{theorem}

We parameterize triangles by their edge lengths and the radii of the 
discs at these vertices. We label the vertices 
$i=0,1,2$ and we denote the length of the
edge opposite vertex $i$ by $x_i$. The radius of the disc at vertex $i$
is $r_i$. The excess is then written as $E(x_0,x_1,x_2,r_0,r_1,r_2)$, 
and the localizing potential as 
$F(x_0,x_1,x_2,r_0,r_1,r_2)$.
$F$ will be the sum over the three
vertices of the triangle of a {\it vertex localizing  potential}
plus the sum over the three edges of an {\it edge localizing potential}.
The vertex potential is based on the constraint that the 
sum of the angles around a vertex is $2 \pi$.
The edge potential is based on a constraint involving the signed distance 
from the edge to the center of the circle which circumscribes the triangle.

Consider the vertex potential $v$ for vertex $0$.
It depends on $x_0,x_1,x_2$ only through the angle of the 
triangle at vertex $0$, which we denote by $\phi_0$. 
So we write the vertex potential as $v(\phi_0,r_0,r_1,r_2)$. 
The vertex potentials for the other two vertices are 
$v(\phi_1,r_1,r_0,r_2)$ for vertex $1$ and 
$v(\phi_2,r_2,r_0,r_1)$ for vertex $2$.
We will always take $v(\phi_0,r_0,r_1,r_2)$ to be symmetric 
under the interchange of its last two arguments.

Now consider the edge localizing potential for the edge opposite vertex $0$.
It will be a function of the edge lengths $x_0,x_1,x_2$ and the 
radii of the discs at the endpoints of the edge, i.e., $r_1,r_2$. 
We write it as $e(x_0,x_1,x_2,r_1,r_2)$. It is symmetric 
under the simultaneous interchange of $r_1$ with $r_2$ and $x_1$ with $x_2$.  
The total localizing potential for our triangle is 
\bea
F(x_0,x_1,x_2,r_0,r_1,r_2)&=
v(\phi_0,r_0,r_1,r_2)+v(\phi_1,r_1,r_0,r_2)+v(\phi_2,r_2,r_0,r_1) \nonumber \\
&+e(x_0,x_1,x_2,r_1,r_2)+e(x_1,x_0,x_2,r_0,r_2)+e(x_2,x_0,x_1,r_0,r_1)
\eea

Note that each edge in the packing has two edge localizing 
potentials associated with it, while each vertex has $n$ vertex localizing 
potentials associated with it where $n$ is the number of triangles in 
the packing that contain the vertex. Thus the sum of $F(T)$ over all triangles
can be written as the sum of the following two sums. 
The first is the sum over edges of the sum of the 
two edge localizing potentials associated with the edge.
The second is the
sum over vertices of the sum of the $n$ vertex localizing potentials 
associated with the vertex. ($n$ is vertex dependent.)
Thus to prove \reff{global_lp_cond}, it suffices to prove the following 
two conditions. First, for every disc center in the packing we have 
\be
\sum_{i=1}^n v(\theta_i,r,r_i,r_{i+1}) \ge 0 
\label{v_lp_cond}
\ee
where $n \ge 3$ is the number of triangles that have the disc center as a 
vertex, $r$ is the radius of the disc, and 
$r_1,r_2, \cdots,r_n$ are the radii of the discs that surround it.
These are ordered in the natural way, so that 
one triangle has discs of radii $r,r_1,r_2$, the next has radii of 
$r,r_2,r_3$ and so on. $r_{n+1}$ is defined to be $r_1$. 
$\theta_1,\theta_2,\cdots,\theta_n$ are the angles of these triangles at
the vertex. 
Second, for every edge in the packing we have 
\be
e(x_0,x_1,x_2,r_1,r_2)+e(x_0,x^\prime_1,x^\prime_2,r_1,r_2) \ge 0
\label{edge_lp_cond}
\ee
Here the length of the edge is $x_0$, and $r_1,r_2$ are the radii of 
the two discs at its endpoints. $x_1,x_2$ are the lengths of the other
two edges in one triangle and $x^\prime_1,x^\prime_2$ are 
the lengths of the other two edges in the other triangle. 
The edge potential $e(x_0,x_1,x_2,r_1,r_2)$ will depend on 
the signed distance from the ``center'' of the triangle to the edge opposite
vertex $0$. 
The center of the triangle is the point equidistant from the 
three vertices. The signed distance is 
positive when the center lies on the same side of the edge as the 
vertex opposite the edge. 
For two triangles in a Delaunay decomposition that share an edge, there is 
a constraint on the two signed distances to this common edge that 
will be the basis for proving \reff{edge_lp_cond}. See section 
\ref{section_edge_lp}.

We end this section with a discussion of why we use the Delaunay 
triangulation rather than the FM triangulation \cite{FM}
that was used by Heppes in his proof of the optimality of six other 
compact packings. For most triangles the excess is positive. 
As we will see in the next section, it can be negative 
for triangles that are close to some of the 
triangles that appear in the compact packing of figure \ref{fig545}.
The only other triangles with negative excess are relatively flat
triangles, i.e., triangles with a large obtuse angle. 
The FM triangulation has the nice property that 
each disc is covered by the triangles that have a vertex at the 
center of the disc, and so triangles cannot be too flat. However, 
an FM triangulation can still contain triangles which are flat enough that 
their excess is slightly negative. The Delaunay triangulation
can contain triangles that are even flatter and so have an even 
more negative excess. But in the Delaunay triangulation it is possible 
to define an edge localizing potential that takes care of these 
flat triangles. This is explained in section \ref{section_edge_lp}.

\section{The vertex localizing potentials}
\label{section_vertex_lp}

We begin the proof of the theorem 
by considering what the values of the localizing potentials
must be for the triangles that appear in our densest packing. 
There are four triangles that appear in this packing. 
We will refer to the triangle that has two large discs and one small disc
as the alpha triangle. 
It has one side of length $2$ and two of length $1+r$. We denote the 
angle opposite the side of length $2$ by $\alpha$ and the other two angles by 
$\alpha^\prime$. 
For this triangle the excess is negative, 
\be
E_\alpha \approx -0.0022743457
\ee

We will refer to the triangle that has one large disc and two small discs
as the beta triangle. 
It has one side of length $2r$ and two of length $1+r$. The angle opposite the
side of length $2r$ will be called $\beta$, and the other two are 
$\beta^\prime$. 
This triangle also has negative excess. 
\be
E_\beta \approx -0.0017217279 
\ee

The triangle with three small discs will be called the 
small equilateral triangle. Its sides are of length $2r$, and its excess is 
positive.
We denote it by 
\be
E_S \approx 0.0024336170  
\ee

The fourth triangle has three large discs and will be called the 
large equilateral triangle.
Its sides are of length $2$, and its excess is also positive.
We denote it by 
\be
E_L \approx 0.0081887688 
\ee

The localizing potential must be defined so that it equals the excess for
each of the four triangles. The edge localizing potential is zero 
for all four of the triangles. So we have the 
following four conditions. 
\be
v(\alpha,r,1,1) + 2 v(\alpha^\prime,1,r,1) = E_\alpha
\ee
\be
v(\beta,1,r,r) + 2 v(\beta^\prime,r,r,1) = E_\beta
\ee
\be
3 v(\pi/3,r,r,r) = E_S
\ee
\be
3 v(\pi/3,1,1,1) = E_L
\ee
The localizing potential condition \reff{v_lp_cond} for a small disc
requires
\be
2 v(\alpha,r,1,1) + 2 v(\beta^\prime,r,1,r) + v(\pi/3,r,r,r)=0
\label{lpsmall}
\ee
and for a large disc it requires
\be
v(\beta,1,r,r) + 4 v(\alpha^\prime,1,r,1) + 2 v(\pi/3,1,1,1)=0
\label{lplarge}
\ee

We have found six conditions above, but 
these six conditions are not linearly independent.
In the packing shown in figure \ref{fig545}, 
the unit cell has 6 alpha triangles, 
3 beta triangles, 1 small equilateral
triangle and 2 large equilateral triangles. The sum of the excesses 
of the triangles in a unit cell must be zero, so 
\be
6 E_\alpha + 3 E_\beta + E_S + 2E_L =0
\ee
This implies that the above six conditions are equivalent to five linearly 
independent conditions.

We solve these equations by first introducing two parameters $x$ and 
$y$. We set 
\bea
v(\alpha,r,1,1) &=x E_\alpha \nonumber \\
v(\alpha^\prime,1,r,1) &= {1 \over 2} (1-x) E_\alpha \nonumber \\
v(\beta,1,r,r) &=y E_\beta \nonumber \\
v(\beta^\prime,r,r,1) &= { 1 \over 2} (1-y) E_\beta \nonumber \\
v(\pi/3,r,r,r) &= {1 \over 3} E_S \nonumber \\
v(\pi/3,1,1,1) &= {1 \over 3} E_L 
\label{eqmin}
\eea
The conditions \reff{lpsmall} and \reff{lplarge}
are equivalent conditions on $x$ and $y$. They give
\be x 2 E_\alpha - y E_\beta = 2 E_\alpha + {2 \over 3} E_L
= - E_\beta - {1 \over 3} E_S
\ee
We take $x=0$, and then $y$ is determined by the above equation. 
We will explain the motivation for this choice at the end of this 
section.

We let $\phi_0(r_0,r_1,r_2)$ denote the angle $\phi_0$ in the triangle 
with discs of radius $r_0,r_1,r_2$. So it equals 
$\alpha,\alpha^\prime,\beta,\beta^\prime$ or $\pi/3$. 
Then the vertex localizing potential is 
\be 
v(\phi,r_0,r_1,r_2) = v(\phi_0(r_0,r_1,r_2),r_0,r_1,r_2) + 
   m |\phi-\phi_0(r_0,r_1,r_2)|
\label{vvertex_def}
\ee
where $v(\phi_0(r_0,r_1,r_2),r_0,r_1,r_2)$ is given by the above equations.
So equations \reff{eqmin} are minima of the localizing potential, 
and the potential increases with slope $\pm m$ as we move away from 
these minima. 
Note that we use the same slope, $m$, for all the localizing potentials.
We will take $m=0.12$.

Next we prove that the vertex localizing potential we have defined 
satisfies \reff{v_lp_cond}.
We consider a disc $D$ and must show 
\be
\sum_{i=1}^n v(\theta_i,\rho,r_i, r_{i+1}) \ge 0 
\label{vertex_lp_cond}
\ee
Here $n \ge 3$ is the number of triangles with a vertex at the center of $D$.
$\rho$ is the radius of $D$, and 
$r_1,r_2, \cdots,r_n$ are the radii of the discs that surround $D$.
The angles $\theta_1,\theta_2,\cdots,\theta_n$ 
are the angles of these triangles at
$D$. Of course, the sum of the angles around a vertex is $2 \pi$. So 
\be
\sum_{i=1}^n \theta_i = 2 \pi
\ee

For $\rho=r$, let $n_{rrr}$ be the number of $i$ with $r_i=r_{i+1}=r$.
Let $n_{r11}$ be the number of $i$ with $r_i=r_{i+1}=1$.
And let $n_{r1r}$ be the number of $i$ with one of $r_i,r_{i+1}$ equal 
to $r$ and one equal to $1$. (Keep in mind the symmetry 
$v(\phi,r_0,r_1,r_2)=v(\phi,r_0,r_2,r_1)$.) 
Then to prove inequality \reff{vertex_lp_cond} for $\rho=r$ it 
suffices to prove 
\be
n_{rrr} {E_S \over 3} + n_{rr1} { 1-y \over 2} E_\beta 
+ n_{r11} x E_\alpha 
+ m |2 \pi - n_{rrr} {\pi \over 3} - n_{rr1} \beta^\prime - n_{r11} \alpha | 
\ge 0
\label{nconds}
\ee
And to prove the inequality for $\rho=1$ it suffices to show
\be
n_{1rr} y E_\beta + n_{1r1}{1-x \over 2} E_\alpha+n_{111} {E_L \over 3}
+ m |2 \pi- n_{1rr}\beta - n_{1r1} \alpha^\prime - n_{111} {\pi \over 3}| 
\ge 0
\label{ncondl}
\ee
where $n_{1rr}$, $n_{1r1}$ and $n_{111}$ are defined in the obvious way.
Each choice of $n_{rrr}, n_{rr1}, n_{r11}$ gives a lower bound on $m$,
as does each choice of $n_{1rr}, n_{1r1}, n_{111}$.
Note that for $n_{rrr}=1,n_{rr1}=2, n_{r11}=2$, the left side of 
\reff{nconds} is zero, and for $n_{1rr}=1,n_{1r1}=4, n_{111}=2$, 
the left side of \reff{ncondl} is zero. These are the cases which 
occur in the optimal packing. The largest lower bound on $m$ that we find
from the other cases is when $n_{1rr}=4, n_{1r1}=4, n_{111}=0$.
This case implies
\be
m \ge 0.1185912
\ee
We take 
\be
m = 0.12
\ee
The above computation is the motivation for the choice of $x=0$. 
Other choices force $m$ to be larger.

Finally, we make a modification to definition \reff{vvertex_def}.
We define $v(\phi,r_0,r_1,r_2)$ to be given by the above definition
provided the value is less than $0.1$. Otherwise we define it to be $0.1$. 
Proving \reff{v_lp_cond} for this modified function is easy. 
The most negative value of $v(\phi,r_0,r_1,r_2)$ is 
$E_\alpha \approx -0.0022743457$. So if one or more of 
the $v$ in \reff{v_lp_cond} is $0.1$, then \reff{v_lp_cond} is trivially
satisfied. If all the $v$ in \reff{v_lp_cond} are given by \reff{vvertex_def}, 
then the previous proof applies. 

\section{ The edge localizing potentials}
\label{section_edge_lp}

In this section we define the edge localizing potential and show that 
its sum over the triangles is non-negative. The edge localizing 
potential is only nonzero for relatively flat triangles. 
Such triangles can have negative excess, but in the Delaunay triangulation
they will be adjacent to a triangle with positive 
excess. The edge localizing potential exploits this fact.

We will refer to the point equidistant to the three vertices of a 
triangle as the ``center'' of the triangle. Note that it need not lie
inside the triangle. 
The edge potential $e(x_0,x_1,x_2,r_1,r_2)$ will depend on 
the signed distance from the center of the triangle to the edge opposite
vertex $0$. We let $d(x_0,x_1,x_2)$ denote this signed distance. 
We define it to be positive if the 
center lies on the same side of the edge as vertex $0$, and negative if they
lie on opposite sides. We give the formula for this signed distance 
in appendix \ref{appendix_formula}. We note that the signed distance for 
an edge is negative when the angle opposite the edge is obtuse.

Throughout this section we will consider two triangles that 
share an edge. The length of this common edge is $x_0$,
and $r_1,r_2$ are the radii of the two discs centered at its endpoints. 
$x_1,x_2$ are the lengths of the other
two edges in one triangle and $x^\prime_1,x^\prime_2$ are 
the lengths of the other two edges in the other triangle. 

Now consider $d(x_0,x_1,x_2)$ and $d(x_0,x^\prime_1,x^\prime_2)$, the 
signed distances from the centers of the two triangles to their common edge. 
We claim that if the two triangles come from 
a Delaunay decomposition, then 
\be 
d(x_0,x_1,x_2) + d(x_0,x^\prime_1,x^\prime_2) \ge 0 
\label{dist_cond}
\ee
This inequality is obviously not true for arbitrary triangles that share 
an edge. It says that if the signed distance from the center of a 
triangle to an edge is negative, then the signed distance from the 
center of the other triangle that shares this edge to the edge must be 
positive and greater in absolute value. 
Given this inequality we can take any function $f(d,x_0,r_1,r_2)$ 
which is an increasing and odd function of $d$, and let 
\be
e(x_0,x_1,x_2,r_1,r_2)=f(d(x_0,x_1,x_2),x_0,r_1,r_2)
\ee
Then \reff{dist_cond} implies \reff{edge_lp_cond}. 

\begin{figure}[tbh]
\includegraphics{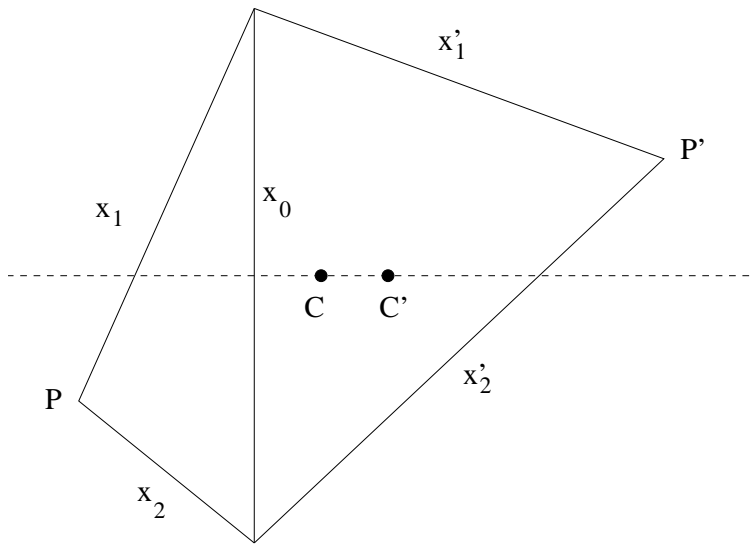}
\caption{$C$ and $C^\prime$ are the centers of the left and right triangles
respectively. In this example the signed distances from $C$ and $C^\prime$ 
to the common edge are negative and positive, respectively.}
\label{fig_dist}
\end{figure}

To prove \reff{dist_cond} we use figure \ref{fig_dist}.
We have drawn the two triangles so their common edge is vertical, 
and we have shown the bisector of this common edge with a dashed line.
The two centers are each equidistant from the endpoints of this common 
edge. So both centers lie on the dashed line.
Inequality \reff{dist_cond} is equivalent to 
$C'$ being to the right of $C$ (or equal to $C$).
By definition the three vertices of a triangle are equidistant to 
the center of the triangle. 
The Delaunay decomposition has the property that no other vertex of a triangle 
is closer to the center than these three vertices. 
In particular, $C$ is closer to $P$ than to $P^\prime$, and 
$C'$ is closer to $P'$ than to $P$. 
It follows that $C'$ is to the right of $C$ (or equal to it).

A fairly simple function $f(d,x_0,r_1,r_2)$ will suffice for our purposes.
Recall that each of $r_1$ and $r_2$ is either $r$ or $1$. 
We define
\be
 f(d,x_0,r,r)=
 \cases{0,&if $x_0<1.8$ \cr
  0.28 \, d,&if $ 1.8 \le x_0 < 2.2$ \cr
  0.4 \, d,&if $ 2.2 \le x_0 $ \cr}
\ee
Next we define
\be
 f(d,x_0,r,1)=
 \cases{0,&if $x_0<2.32$ \cr
  0.06 \, d,&if $ 2.32 \le x_0 $ \cr}
\ee
As always, $f(d,x_0,1,r)=f(d,x_0,r,1)$.
Finally, we let $f(d,x_0,1,1)=0$.

\section{ Local proof the localizing potential works}
\label{section_local}

We now turn to the proof of \reff{local_prob}. 
We will use the computer to prove this inequality for most
triangles. However, for triangles that appear in the densest 
packing shown in figure \ref{fig545}, 
equality holds in \reff{local_prob}.
Thus for triangles close to those that appear in this densest 
packing inequality  \reff{local_prob} will be close to an equality.
So we can only use the computer to prove this inequality for triangles
which do not lie too close to a triangle in the densest packing.
In this section we will prove \reff{local_prob} for the 
triangles which are close to a densest packing triangle. 

There are four triangles that appear in this densest configuration.
In all of these triangles the three discs are tangent to one another. 
As we will see this implies that the excess $E(T)$ 
has a local minimum at each of these triangles. The localizing 
potential $V(T)$ also has a local minimum at these triangles. 
We must show that when we perturb one of these triangles 
the increase in $E(T)$ is greater than the increase in $V(T)$. 

Throughout this section we will obtain bounds on  quantities by 
bounding their partial derivatives with respect to $x_0,x_1$ and $x_2$. 
If $f(x_0,x_1,x_2)$ is a function of the three edge lengths, and 
we have 
\be
 c_i \le {\partial f \over \partial x_i }(x_0,x_1,x_2)  \le d_i
\ee
throughout some neighborhood of $(\bar x_0,\bar x_1,\bar x_2)$, 
then 
\be
\sum_{i=0}^2 c_i \Delta x_i \le 
  f(x_0,x_1,x_2) - f(\bar x_0,\bar x_1,\bar x_2)
\le \sum_{i=0}^2 d_i \Delta x_i 
\ee
where $\Delta x_i=x_i-\bar x_i$. 

We first consider triangles that have discs of radius $r$ at all three 
vertices and which satisfy
\be 
2 r \le x_i \le 2 r + \epsilon, \quad i=0,1,2
\label{tri_set_equia}
\ee
where $\epsilon>0$ will be determined later. 
Since the three discs all have radius $r$, $D$ is independent of the 
$x_i$. Letting 
\be
\Delta E = E(x_0,x_1,x_2,r,r,r)- E(2r,2r,2r,r,r,r),
\ee
then
\be
\Delta E = \Delta(\delta A - D) 
 \ge \delta \sum_{i=0}^2 a^i_{rrr}(\epsilon) \Delta x_i
 = \delta a_{rrr}(\epsilon) \sum_{i=0}^2 \Delta x_i
\ee
where 
\be
a_{rrr}(\epsilon) = \min {\partial A \over \partial x_i}
\ee
with the min over the above set of triangles. 
Note that this min is independent of $i$. 
Inequality \reff{a0_partial} in the appendix gives a lower bound on 
$a_{rrr}(\epsilon)$. 

We must compare the above with the increase in the localizing potential
$V(T)$. Let $\Delta \phi_i$ denote the change in the angles $\phi_i$ 
corresponding to changing each $x_j$ by $\Delta x_j$.
Then the change in $V(T)$ is bounded by 
\be
\Delta V(T) \le m \sum_{i=0}^2 |\Delta \phi_i|
\le m \sum_{j=0}^2 b^j_{rrr}(\epsilon) \Delta x_j
\ee
where 
\be
b^j_{rrr}(\epsilon) = \max \sum_{i=0}^2 |{\partial \phi_i \over \partial x_j}|
\ee
and the max is over the set of triangles. 
$b^i_{rrr}(\epsilon)$ is independent of $i$, and we denote it by 
$b_{rrr}(\epsilon)$.
This proves that the excess is greater than the localizing potential 
for all triangles satisfying \reff{tri_set_equia} provided 
$m \le m_{rrr}(\epsilon)$ where 
\be
m_{rrr}(\epsilon) = \delta \, {a_{rrr}(\epsilon) \over b_{rrr}(\epsilon)} 
\label{mbounda}
\ee

For triangles that have discs of radius $1$ at all three 
vertices the estimates are very similar.
We consider the set of triangles which satisfy
\be 
2 \le x_i \le 2 + \epsilon, \quad i=0,1,2
\label{tri_set_equib}
\ee
The constraint we obtain on $m$ is now $m \le m_{111}(\epsilon)$ where 
\be
 m_{111}(\epsilon) = \delta \, {a_{111}(\epsilon) \over b_{111}(\epsilon)}
\label{mboundb}
\ee
$a_{111}(\epsilon)$ and $b_{111}(\epsilon)$ are defined by analogy to 
$a_{rrr}(\epsilon)$ and $b_{rrr}(\epsilon)$.

Next we consider triangles with one disc of radius $1$ and two of 
radius $r$. We take the large disc to be at vertex $0$ and consider the 
set of triangles given by 
\be 
2r \le x_0 \le 2r + \epsilon, \quad
1+r \le x_i \le 1+r + \epsilon, \quad i=0,1,2
\label{tri_set_equic}
\ee
In this case $D$ is not independent of the $x_i$. 
Using the constraint $\phi_0+\phi_1+\phi_2=\pi$, we can write $D$ as   
\be
D={1 \over 2} (\phi_0+ r^2 \phi_1 + r^2 \phi_2) 
={\pi \over 2} r^2  + {1 - r^2\over 2} \phi_0
\ee
Thus 
\be
\Delta E = \Delta(\delta A - D) 
 \ge \sum_{i=0}^2 \delta a^i_{1rr}(\epsilon) \Delta x_i 
  - {1 - r^2\over 2} \Delta \phi_0
\ee
where $a^i_{1rr}$ is the min of ${\partial A \over \partial x_i}$
over the set of triangles.

Let 
\bea
c^i_{1rr}(\epsilon) = \min {\partial \phi_0 \over \partial x_i} \nonumber \\
d^i_{1rr}(\epsilon) = \max {\partial \phi_0 \over \partial x_i} \nonumber \\
\eea
where the $\min$ and $\max$ are over the triangles satisfying 
\reff{tri_set_equic}. Note that there are no absolute values on the 
partial derivatives in the above. The calculations in the appendix
show that on the set of triangles ${\partial \phi_0 \over \partial x_i}$ 
is positive for $i=0$ and negative for $i=1,2$. 
We now have 
\be
\Delta E \ge \sum_{i=0}^2 \left[\delta a^i_{1rr}(\epsilon) 
  - {1 - r^2\over 2} d^i_{1rr}(\epsilon) \right] \Delta x_i
\ee
We bound the increase in the localizing potential $V(T)$ as before. 
The constraints we obtain on $m$ are $m \le m^i_{1rr}(\epsilon)$ with 
\be
m^i_{1rr}(\epsilon) = 
{\delta a^i_{1rr}(\epsilon) - {1 - r^2\over 2} d^i_{1rr}(\epsilon) \over
b^i_{1rr}(\epsilon)}
\label{mboundc}
\ee

Finally we consider triangles with one disc of radius $r$ and two of 
radius $1$. We take the small disc to be at vertex $0$ and consider the 
set of triangles given by 
\be 
2 \le x_0 \le 2 + \epsilon, \quad
1+r \le x_i \le 1+r + \epsilon, \quad i=0,1,2
\label{tri_set_equid}
\ee
We now have 
\be
D={1 \over 2} (r^2 \phi_0+ \phi_1 + \phi_2) 
={\pi \over 2} r^2  + {r^2 -1 \over 2} \phi_0
\ee
Thus
\be
\Delta E \ge \sum_{i=0}^2 \delta a^i_{r11}(\epsilon) \Delta x_i 
  + {1 - r^2\over 2} \Delta \phi_0
\ge \sum_{i=0}^2 \left[ \delta a^i_{r11}(\epsilon) 
  + {1 - r^2\over 2} c^i_{r11}(\epsilon) \right] \Delta x_i
\ee
and so we obtain the constraints $m \le m^i_{r11}(\epsilon)$ with 
\be
m^i_{r11}(\epsilon) 
= {\delta a^i_{r11}(\epsilon) + {1 - r^2\over 2} c^i_{r11}(\epsilon) \over
b^i_{r11}(\epsilon)}
\label{mboundd}
\ee

Using the results of appendix \ref{appendix_near}, we can evaluate 
these bounds on $m$. We will take $\epsilon=0.001$. 
For comparison we also show the values for $\epsilon=0$. 
We find 
\bea 
& m_{rrr}(0)=0.135463, \quad m_{rrr}(0.001)=0.134576 \nonumber \\
& m_{111}(0)=0.455814, \quad m_{111}(0.001)=0.454183 \nonumber \\
& m^0_{1rr}(0)=0.232960, \quad m^0_{1rr}(0.001)=0.231100 \nonumber \\
& m^i_{1rr}(0)=0.179205, \quad m^i_{1rr}(0.001)=0.178067 \nonumber \\
& m^0_{r11}(0)=0.264015, \quad m^0_{r11}(0.001)=0.262815  \nonumber \\
& m^i_{r11}(0)=0.308628, \quad m^i_{r11}(0.001)=0.306788 \nonumber \\
\eea
where $i$ is $1$ or $2$. 
Recall that we take $m=0.12$, so these bounds are all easily met. 
In fact, with $m=0.12$ this local proof works for $\epsilon$ as large 
as $0.018$. 

\section{ Global  proof the localizing potential works}
\label{section_global}

We now prove inequality \reff{local_prob} for the rest of the triangles,
i.e., all the triangles that are not close 
to those that appear in the densest configuration in figure \ref{fig545}.
We use interval arithmetic for this proof.
This is necessary for two reasons. First, it is needed (as in all 
computer-assisted proofs involving real numbers) to deal with 
the fact that computer calculations involving real numbers are not
exact. Second, even if the computer could perform exact calculations
we could not prove \reff{local_prob} by checking one triangle at a 
time since there are infinitely many triangles. Instead we must
work with sets of triangles. Given intervals for the three edge
lengths $x_0,x_1,x_2$ of the triangle, we consider the set of triangles
whose edge lengths belong to the respective intervals. It is 
straightforward to use interval arithmetic to then compute intervals 
for quantities such as the angles $\phi_0,\phi_1,\phi_2$. Interval
arithmetic is designed so that the interval computed for 
$\phi_i$ means that for a triangle whose edge lengths 
belong to their respective intervals, the value of $\phi_i$ must belong
to its interval. We note that the code for computing intervals 
for quantitites such as $\phi_i$ is the same as the code for computing
these quantities when we have a single real value for each edge length.
We need only change the data type of the variables involved from 
``double'' to ``interval'' and define versions of the basic arithmetic
operations (addition, multiplication, inverse cosine, etc.) for 
interval variables. The localizing potential involves real numbers 
that must be represented by intervals, and so it too is handled using
interval arithmetic. 

We now fix choices of $r_0,r_1,r_2$ and ask what triangles must be 
considered in the proof of \reff{local_prob}.
Since discs cannot overlap, we have the following lower bounds on 
the length of the edges of the triangle.
\bea
x_0 &\ge r_1+r_2 \nonumber \\
x_1 &\ge r_0+r_2 \nonumber \\
x_2 &\ge r_0+r_1 \nonumber \\
\eea
Recall that the center of a triangle is the point equidistant to the 
three vertices, and the circumradius is the radius of the circle 
centered at this point which contains all three vertices.
In the Delaunay decomposition no disc has a center closer to the triangle's
center than the discs at its vertices. So if the circumradius is greater
than $1+r$, then we can add a disc of radius $r$ to the packing 
by putting its center at the center of the triangle. If we assume
that the packing is saturated, then this implies that every triangle 
has a circumradius of at most $1+r$. The circumradius is greater than half
the length of any side of the triangle. So we have 
\be
x_i \le 2(1+r), \quad i=0,1,2
\ee

At each stage of the computer proof we have a list of 
parallelepipeds $I_0 \times I_1 \times I_2$. Such a parallelepiped 
represents the set of triangles with $x_i \in I_i$. 
At each stage of the computer proof the union of the parallelepipeds 
in the list is the set of triangles for which we must still prove 
\reff{local_prob}. Initially, there is one parallelpiped in the list
with 
\be 
I_0=[r_1+r_2,2(1+r)], \quad I_1=[r_0+r_2,2(1+r)], \quad I_2=[r_0+r_1,2(1+r)]
\ee
At each stage we take a parallelepiped from the list and split it 
into two parallelepipeds. We compute the 
interval for $E(T) - F(T)$ for each of the two. If it contains 
no negative values for a particular parallelepiped, then we know 
that \reff{local_prob} is true for all triangles in this parallelepiped, and 
we can discard it. 
If the interval contains some negative values, we 
must add this parallelpiped to the list.
Note that if the parallelepiped is contained 
in one of the parallelepipeds \reff{tri_set_equia}, \reff{tri_set_equib}, 
\reff{tri_set_equic}, or \reff{tri_set_equid}, then we know $E(T)-F(T) \ge 0$
on this parallelepiped, and so it can be discarded.  
In appendix \ref{appendix_formula} we show that if the packing is 
saturated then in the Delaunay triangulation the area of every triangle 
is at least $2 r^3 /(1+r)$. So for a parallelepiped, if the interval 
we compute for the area of the triangles is entirely greater than this 
lower bound, then we can discard the parallelepiped.
The process ends when the list of parallelepipeds for which 
\reff{local_prob} has not been proved is empty.

We compute our interval bounds in the simplest, crudest way. 
For example, we compute the interval bound on the angle $\phi_i$ by 
simply using interval arithmetic in the equation for $\phi_i$ in 
terms of the edge lengths, eq. \reff{phieq}. One could compute 
a better interval bound by explicitly maximizing and minimizing
$\phi_i$ over the parallelepiped by considering the derivatives 
of $\phi_i$ with respect to the edge lengths $x_j$.
Because of our rather crude bounds, the parallelepipeds must be
split a large number of times. The computer considers a total of about 
26 million parallelepipeds. This takes about 30 minutes on a laptop.

\begin{appendix}

\section{Formulae}
\label{appendix_formula}

The vertices are labelled 0,1,2. 
The length of the edge opposite vertex $i$ is $x_i$. The radius 
of the disc at vertex $i$ is $r_i$. Each $r_i$ can only be $r$ or $1$. 
The angle at vertex $i$ is $\phi_i$. It is given by 
\be
 \cos(\phi_i)= (x_j^2 + x_k^2 - x_i^2)/(2 x_j x_k)
\label{phieq}
\ee
where $\{i,j,k\}=\{0,1,2\}$.
The area of the triangle is 
\be
A={1 \over 4} \sqrt{2 x^2_0 x^2_1 + 2 x^2_0 x^2_2 + 2 x^2_1 x^2_2 
- x_0^4 - x_1^4 - x_2^4}
\ee

Recall that by the center of a triangle we mean the center of the 
circle which contains the three vertices of the triangle. 
We will derive a formula for the signed distance from the center 
to each edge and a formula for the circumradius. 
We take the vertices of the triangle to be 
\be 
P_0 = (0,0), \quad P_1 = (x_0,0) \quad P_2 = (a,b) 
\ee
where $a^2+b^2=x_2^2$. 
Let $P$ be the center of the triangle. 
The center must lie on the line $x=x_0/2$, so $P=(x_0/2,d_0)$, where
$d_0$ is the signed distance from the center to the edge opposite vertex
$0$.  

The circumradius of the triangle, $R$, satisfies
\be
 R=||P-P_0||=||P-P_1||
\ee
This implies
\be 
R^2={1 \over 4} x_0^2 + d_0^2 = (a-x_0/2)^2 + (b-d_0)^2
\label{circ_eq}
\ee
Using $a^2 + b^2 = x_2^2$, this gives
\be
a x_0 + 2 b d_0 = x_2^2
\ee
Using
\be
a=x_2 \cos(\phi_1) = x_2 {x_0^2+x_2^2 - x_1^2 \over 2 x_0 x_2 }
\ee
and 
\be 
{1 \over 2} b x_0 = A
\ee
we find 
\be
d_0={x_0(x_1^2+x_2^2-x_0^2) \over 8 A }
\ee

Using \reff{circ_eq}, the above formula for $d_0$ yields a formula
for the circumradius.
\be
R= {x_0 x_1 x_2 \over 4 A}
\ee
For a saturated packing we can assume 
\be
R \le 1+r
\ee
Otherwise we could add a disc of radius $r$ with center at the center of 
the triangle. Using the above formula for $R$ this implies a lower bound on 
the area $A$. 
\be
A \ge {x_0 x_1 x_2 \over 4 (1+r) }  
\ee
Using the trivial bound $x_i \ge 2 r$, this implies  
\be
A \ge {2 r^3 \over 1+r} 
\label{area_bound}
\ee

\section{Calculations for triangles with excess near zero}
\label{appendix_near}

We consider triangles with $x_0 \in [x,x+\epsilon]$, 
$x_1,x_2 \in [y,y+\epsilon]$.

We have 
\be
16 A^2 = 2 x^2_0 x^2_1 + 2 x^2_0 x^2_2 + 2 x^2_1 x^2_2 
- x_0^4 - x_1^4 - x_2^4
\ee
Given an $i=0,1$ or $2$, we let $j$ and $k$ be the other two integers in 
$\{0,1,2\}$. 
Then 
\be
8 A {\partial A \over \partial x_i} = 
x_i(x_j^2+x_k^2-x_i^2)
\ee
We assume that $\epsilon$ is small enough that the set of triangles
we are considering only contains acute triangles. So 
$x_j^2+x_k^2-x_i^2>0$. Then ${\partial A \over \partial x_i} \ge 0$ 
for each $i$, and so $A$ is increasing in each $x_i$. So $A$ is bounded 
above by its value for the triangle with sides of length 
$x_0=x+\epsilon$, $x_1=x_2=y+\epsilon$.
This area is 
\be
A_{max} = {1 \over 4} (x+\epsilon) \sqrt{4(y+\epsilon)^2 - (x+\epsilon)^2}
\ee
Thus we obtain the lower bounds 
\be
{\partial A \over \partial x_0} 
\ge { x (2 y^2-(x+\epsilon)^2) \over 8 A_{max} }
\label{a0_partial}
\ee
and for $i=1,2$
\be
{\partial A \over \partial x_i}
\ge { y(x^2+y^2-(y+\epsilon)^2) \over 8 A_{max} }
\label{a1_partial}
\ee

From \reff{phieq} we have 
\be
{\partial \phi_0 \over \partial x_0} = (\sin(\phi_0))^{-1} \, 
{x_0 \over x_1 x_2}
\ee
\be
{\partial \phi_0 \over \partial x_1} = - (\sin(\phi_0))^{-1} \,
{x_1^2 + x_0^2 - x_2^2 \over 2 x_1^2 x_2}
\ee
\be
{\partial \phi_0 \over \partial x_2} = - (\sin(\phi_0))^{-1} \,
{x_2^2 + x_0^2 - x_1^2 \over 2 x_2^2 x_1}
\ee
Note that ${\partial \phi_0 \over \partial x_0}$ is always positive 
and the other two partials are negative for acute triangles. 
So $\sin(\phi_0)$ is bounded below by its
value at $x_0=x$, $x_1=x_2=y+\epsilon$.
Let $S_0$ be the value of  $\sin(\phi_0)$ at this point. 
And $\sin(\phi_0)$ is bounded above by its
value at $x_0=x+\epsilon$, $x_1=x_2=y$.
Let $S_0^\prime$ be the value of  $\sin(\phi_0)$ at this point. 
Then we have 
\be
 {x \over (y+\epsilon)^2 S_0^\prime} \le
{\partial \phi_0 \over \partial x_0} \le {x+\epsilon \over y^2 S_0}
\ee
\be
- {(y+\epsilon)^2  + (x+\epsilon)^2  - y^2 \over 2 y^3 S_0} \le 
{\partial \phi_0 \over \partial x_i} \le 
- {y^2  + x^2  - (y+\epsilon)^2 \over 2 (y+\epsilon)^3 S_0^\prime} 
\ee
for $i=1,2$.

Similarly
\be
{\partial \phi_1 \over \partial x_1} = (\sin(\phi_1))^{-1} \, 
{x_1 \over x_0 x_2}
\ee
\be
{\partial \phi_1 \over \partial x_0} = - (\sin(\phi_1))^{-1} \,
{x_0^2 + x_1^2 - x_2^2 \over 2 x_0^2 x_2}
\ee
\be
{\partial \phi_1 \over \partial x_2} = - (\sin(\phi_1))^{-1} \,
{x_2^2 + x_1^2 - x_0^2 \over 2 x_2^2 x_0}
\ee
For these derivatives we only need bounds on their absolute values.
\be
|{\partial \phi_1 \over \partial x_1}| \le  
{(y+\epsilon) \over x y S_1}
\ee
\be
|{\partial \phi_1 \over \partial x_0}| \le 
{(x+\epsilon)^2 + (y+\epsilon)^2 - y^2 \over 2 x^2 y S_1}
\ee
\be
|{\partial \phi_1 \over \partial x_2}| \le 
{2(y+\epsilon)^2  - x^2 \over 2 y^2 x S_1}
\ee
with $S_1$ equal to the value of $\sin(\phi_1)$ when 
$x_1=y$, $x_0=x+\epsilon$ and $x_2=y+\epsilon$. 

\section{ Equation for r}
\label{appen_eqr}

In this appendix we show that the radius $r$ for the compact packing
is a root of an eighth degree polynomial. The angles satisfy 
\be
2 \alpha - \beta = 120
\ee
So 
\be
\cos(\alpha-\beta/2)=\cos(60)=1/2
\ee
Hence
\be
\cos(\alpha) \cos(\beta/2)+\sin(\alpha) \sin(\beta/2)=1/2
\ee
This implies 
\be
\sin^2(\alpha) \sin^2(\beta/2) = (1/2 - \cos(\alpha) \cos(\beta/2))^2
\ee
which is equivalent to 
\be
1-\cos^2(\alpha)^2 - \cos^2(\beta/2) = {1 \over 4} - \cos(\alpha) \cos(\beta/2)
\ee
We have
\be
\cos(\alpha)={2(1+r)^2-2^2 \over 2 (1+r)^2} = {r^2+2r-1 \over (1+r)^2}
\ee
and 
\be
\sin(\beta/2)={r \over 1+r}
\ee
\be
\cos(\beta/2)={\sqrt{2r+1} \over 1+r}
\ee
Substituting these into the above equation, some algebra gives
\be
r^8 - 8r^7 - 44 r^6 -232 r^5 -482 r^4
-24 r^3+388 r^2 -120 r + 9 = 0 
\ee

By considering the packing itself, the conditions that various
pairs of discs are tangent can be used to show that $r$ also 
satisfies 
\be
(7+4\sqrt{3}) r^4 + (20+12 \sqrt{3}) r^3 + 
(6+4\sqrt{3}) r^2+(-20-4\sqrt{3})r+3=0
\ee

To compute the density of the packing, consider the unit cell 
shown in figure \ref{fig545}. The quadrilateral has sides of 
length 
\be
l=1+r+\sqrt{3} r + \sqrt{1+2r}
\ee
and angles of $\pi/2$ and $2\pi/3$. So its area is $\sqrt{3} l^2 /2$. 
It covers a total of $3$ large discs and $3$ small discs. So 
the density is 
\be
\delta= {2 \sqrt{3} \, \pi (1+r^2) \over l^2} \approx 0.911627478    
\ee

\end{appendix}

\section*{Acknowledgements}
The author would like to thank Aladar Heppes for helpful correspondence. 
This work was supported by the National Science Foundation (DMS-0201566).

\end{document}